\documentclass{article}

\newtheorem{fed}{\textbf{Definition}}[section]

\newtheorem{lemma}[fed]{\textbf{Lemma}}

\usepackage{amssymb,bbm,graphicx,epsfig,psfrag,epic,eepic,latexsym}
\usepackage{amsmath}
\usepackage{mathrsfs}
\usepackage[dvips]{color}

\begin{document}
\title{Delayed Rabinowitz Floer Homology}
\author{Urs Frauenfelder}
\maketitle

\begin{abstract}
In this article we study Rabinowitz Floer Homology for several interaction
particles. In general Rabinowitz action functional is invariant under simultaneous
time translation for all particles but not invariant if the times of each particle
are translated individually. The delayed Rabinowitz action functional is invariant
under individual time translation for each particle. Although its critical point
equation looks like a Hamiltonian delay equation it is actually an ODE in disguise and
nothing else than the critical point equation of the undelayed Rabinowitz action functional. We show that we can even interpolate between the two action functionals
without changing the critical points and their actions. Moreover, for each of these
interpolating action functionals we have compactness for gradient flow lines under
a suitable restricted contact type assumption. 
\end{abstract}

\section{Introduction}

Rabinowitz Floer homology is the semi-infinite dimensional Morse homology in the
sense of Floer associated to Rabinowitz action functional
\cite{cieliebak-frauenfelder1}. Rabinowitz action functional is the Lagrange multiplier
functional for minus the area functional to the constraint given by the mean value
of a Hamiltonian. It detects periodic orbits of this Hamiltonian for fixed energy
but arbitrary period. The periodic orbits are parametrized. To get
unparametrized periodic orbits one has to mod out the circle action given by
reparametrisation of the free loop space. 
\\ \\
In this article we are interested in several interacting particles. The phase space
of each particle is modelled by a symplectic manifold $M_i$ so that the total phase
space is the product
$$M=\bigoplus_{i=1}^m M_i,$$
when $m$ denotes the number of particles. The Hamiltonian is of the form
$$H=f(H_1,\ldots,H_m)$$
where the $H_i$ are Hamiltonians on $M_i$ and $f \colon \mathbb{R}^m \to \mathbb{R}$
is a smooth function. A motivating example to consider for this kind of set-up are
toric domains. There $M_i=\mathbb{C}$ and $H_i(z)=\pi|z|^2$ is the moment map
for the standard circle action on $\mathbb{C}$ for $1 \leq i \leq m$. It is worth
mentioning in this context that in recent times quite some Hamiltonian systems
of important physical origin where interpreted in terms of toric geometry as
concave toric domains. The pioneering work in this context is the interpretation
of the billiard system on a round table as a concave toric domain by 
Ramos \cite{ramos}. Mohebbi interpreted the rotating Kepler problem below
the first critical value as a concave toric domain \cite{mohebbi}. In 
\cite{frauenfelder} the same was proved for the Stark problem and in view of the proof of the
Dullin-Montgomery conjecture by Pinzari \cite{pinzari} this holds true as well for the Euler problem.
\\ \\
The free loop space of a product of $m$ symplectic manifolds is invariant under
the action of the $m$-dimensional torus $T^m$ which acts by reparametrizing the
loop in each component individually. The periodic orbits of a Hamiltonian as above
are invariant under the $T^m$-action as well. However, Rabinowitz action functional
$\mathcal{A}$ is in general \emph{not} invariant under the torus action, but just under its diagonal circle action, unless the function $f$ is linear. 
\\ \\
To remedy this unpleasant situation we consider in this note a deformation of
Rabinowitz action functional $\mathcal{A}_r$ for $r \in [0,1]$, where
$\mathcal{A}_0=\mathcal{A}$ is the original Rabinowitz action functional while
the action functional $\mathcal{A}_1$ is invariant under the full $T^m$-action. 
The action functional is a Lagrange multiplier version of action functionals
considered in \cite{albers-frauenfelder-schlenk} to study Hamiltonian delay equations.
We therefore refer to the action functional $\mathcal{A}_1$ as the \emph{delayed Rabinowitz action functional}.
Its critical point equation looks like a delay equation. But in this paper we
analyse this equation further and the analysis reveals that it is actually an ODE and nothing else than the periodic orbit equation. In fact the critical points of
$\mathcal{A}_0$ and $\mathcal{A}_1$ coincide. Even more is true. We prove the following
\emph{non-bifurcation theorem}.
\\ \\
\textbf{Theorem\,A: } \emph{The critical set $\mathrm{crit}(\mathcal{A}_r)$ is constant, i.e., independent of $r$, and the restriction of $\mathcal{A}_r$ to it as well.} 
\\ \\
The second result is a compactness result for gradient flow lines of
the functionals $\mathcal{A}_r$. Even for the nondelayed Rabinowitz action
functional $\mathcal{A}_0$ this requires some assumption on the hypersurface
in the symplectic manifold. In \cite{cieliebak-frauenfelder1} compactness was
established under a restricted contact assumption. In Section~\ref{delfund} 
we introduce a restricted contact type assumption for the case of several particles
modelled on a product symplectic manifold. Under this assumption compactness of 
gradient flow lines holds true, namely
\\ \\
\textbf{Theorem\,B: } \emph{Under the restricted contact type assumption,
suppose that $w_\nu$ is a sequence of gradient flow lines for
$\mathcal{A}_{r_\nu}$ where $\nu \in \mathbb{N}$ and $r_\nu \in [0,1]$ such that
there exists an interval $[a,b] \subset \mathbb{R}$ with the property that
$$\mathcal{A}_{r_\nu}(w_\nu(s)) \in [a,b], \quad \nu \in \mathbb{N},\,\,s \in \mathbb{R}.$$
Then $w_\nu$ has a subsequence which converges to a gradient flow line in
the $C^\infty_{\mathrm{loc}}$-topology.}
\\ \\
From Theorem\,A and Theorem\,B it follows that the Rabinowitz Floer homology
of the delayed Rabinowitz action functional is well defined and canonically isomorphic
to the usual Rabinowitz Floer homology. Moreover, since the action according to Theorem\,A stays constant along the critical set the two spectral numbers of the delayed
Rabinowitz Floer homology coincide with the ones of the usual Rabinowitz Floer homology. 
\\ \\
The main interest of the author on this result is that the delayed Rabinowitz 
action functional is invariant under the torus action obtained by changing time for
each particle individually, which is in general not the case for the undelayed
Rabinowitz action functional. This in particular gives us the possibility to
define as well \emph{Tate Rabinowitz Floer homology} for this torus action by using the
delayed Rabinowitz action functional. Spectral numbers for Tate Rabinowitz Floer homology the author is currently studying with Cieliebak \cite{cieliebak-frauenfelder2}
for several harmonic oscillators. The case of several harmonic oscillators corresponds
to the case where the function $f$ is linear so that the delayed Rabinowitz action
functional coincides with the undelayed one. Spectral numbers in Tate Rabinowitz
Floer homology show fascinating connections to the quantum spectrum and therefore play
an important role in the question if there is a homological approach to Gutzwiller's intriguing trace formula \cite{gutzwiller}. In particular, it would be interesting
to address this question for toric domains in view of the close connection referred to above of important Hamiltonian systems arising in atomic physics with toric geometry. 
On the other hand the delayed Rabinowitz action functional studied in this note
can be delayed further to get actual Hamiltonian delay equations as critical points. Such Hamiltonian
delay equation for example show up in the study of Helium for mean interactions
of the electrons \cite{cieliebak-frauenfelder-volkov, frauenfelder0}.
\\ \\
\emph{Acknowledgements: } The author acknowledges partial support by DFG
grant FR 2637/2-2.

\section{The delayed Rabinowitz action functional}

In this section we define the delayed Rabinowitz action functional after having
recalled the undelayed one, show that it is invariant under a torus action 
and proof Theorem\,A from the introduction. 
\\ \\
Suppose that $(M_i,\omega_i=d\lambda_i)$ for $1 \leq i \leq m$ is a finite collection of exact symplectic
manifolds and  $H_i \colon M_i \to \mathbb{R}$ as well as $f \colon \mathbb{R}^m \to
\mathbb{R}$ are smooth functions. On the product manifold
$$M:=\bigoplus_{i=1}^m M_i$$
we consider the smooth function
$$H_f \colon M \to \mathbb{R}, \quad z=(z_1,\ldots,z_m) \mapsto
f\Big(H_1(z_1),\ldots,H_m(z_m)\Big).$$
The product manifold $M$ is itself an exact symplectic manifold with symplectic form
$$\omega=\oplus_{i=1}^m \omega_i \in \Omega^2(M)$$
and primitive
$$\lambda=\oplus_{i=1}^m \lambda_i \in \Omega^1(M).$$
Abbreviate the free loop space of $M$ by
$$\mathcal{L}:=\mathcal{L}_M:=C^\infty(S^1,M),$$
where $S^1=\mathbb{R}/\mathbb{Z}$ is the circle. Rabinowitz action functional
$$\mathcal{A}_0 \colon \mathcal{L} \times \mathbb{R} \to \mathbb{R}$$
at $(v,\tau) \in \mathcal{L} \times \mathbb{R}$
is given by
$$\mathcal{A}_0(v,\tau)=-\int_{S^1}v^* \lambda+\tau \int_0^1 H_f(v(t))dt.$$
We denote by $X_{H_f}$ the Hamiltonian vector field of $H_f$ implicitly defined by
$$dH_f=\omega(\cdot,X_{H_f}).$$
With this notation the critical points of $\mathcal{A}_0$ are solutions of the problem
$$
\left.\begin{array}{cc}
\partial_t v(t)=\tau X_{H_f}(v(t)), & t \in S^1\\
\int_0^1 H_f(v(t))dt=0. & 
\end{array}\right\}
$$
By preservation of energy this is equivalent to the problem
\begin{equation}\label{crit}
\left.\begin{array}{cc}
\partial_t v(t)=\tau X_{H_f}(v(t)), & t \in S^1\\
H_f(v(t))=0 & t \in S^1.
\end{array}\right\},
\end{equation}
i.e., the critical points of $\mathcal{A}_0$ are periodic orbits of
the Hamiltonian vector field $X_{H_f}$ of period $\tau$ of energy zero. The period
$\tau$ is allowed to be zero, which means that the orbit is constant, or negative, which
means that the orbit is traversed backwards. 
\\ \\
Before writing down the delayed Rabinowitz action functional we first rewrite the usual
Rabinowitz action functional a bit. If we write a loop $v \in \mathcal{L}$ into components
$$v=(v_1,\ldots, v_m)$$
where
$$v_i \in \mathcal{L}_{M_i}, \quad 1 \leq i \leq m$$
the Rabinowitz action functional can be written as
$$\mathcal{A}_0(v,\tau)=-\int_{S^1}v^* \lambda+\tau \int_0^1 f\Big(H_1(v_1),\ldots,H_m(v_m)\Big)dt.$$
The \emph{delayed Rabinowitz action functional} 
$$\mathcal{A}_1 \colon \mathcal{L} \times \mathbb{R} \to \mathbb{R}$$
is obtained from the usual Rabinowitz action functional by interchanging the order of the
integral and the function $f$, namely
$$\mathcal{A}_1(v,\tau)=-\int_{S^1} v^* \lambda+\tau f\Bigg(\int_0^1 H_1(v_1)dt,
\ldots \int_0^1 H_m(v_m)dt\Bigg).$$
In order to write this a bit more compactly we introduce the following notation.
We define
$$H \colon M \to \mathbb{R}^m, \quad z=(z_1,\ldots,z_m) \mapsto
\Big(H_1(z_1), \ldots, H_m(z_m)\Big),$$
so that
$$H_f=f \circ H.$$
For any smooth function $G \colon M \to \mathbb{R}^k$ for $k \in \mathbb{N}$
we define the averaged function
$$\overline{G} \colon \mathcal{L} \to \mathbb{R}^k, \quad
v \mapsto \int_0^1 G(v(t))dt.$$
Using these abbreviations we have
$$\mathcal{A}_0(v,\tau)=-\int_{S^1}v^*\lambda+\tau \overline{f \circ H}(v), \quad
\mathcal{A}_1(v,\tau)=-\int_{S^1} v^*\lambda+\tau f \circ \overline{H}(v),$$
i.e., for the usual Rabinowitz action function the function $f$ is applied before
averaging while for the delayed one after averaging. 
\\ \\
The delayed Rabinowitz action functional is invariant under a torus action which 
in general is not true for the undelayed one. On each loop space
$\mathcal{L}_{M_i}$ we have a circle action 
$$S^1 \times \mathcal{L}_{M_i} \to \mathcal{L}_{M_i}$$
given by reparametrisation. Namely
if $v_i \in \mathcal{L}_{M_i}$ and $r \in S^1$ we set
$$r_* v(t)=v(t+r), \quad r \in S^1.$$
This circle action gives rise to a product action of the $m$-dimensional torus
$$T^m:=\underbrace{S^1 \times  \ldots \times S^1}_{m\,\,\mathrm{times}}$$
on $\mathcal{L}_M$ given for $r=(r_1,\ldots,r_m) \in T^m$ and
$v=(v_1,\ldots,v_m) \in \mathcal{L}_M$ by
$$r_*v=\Big((r_1)_*v_1,\ldots,(r_m)_*v_m\Big).$$
We extend this action trivially to $\mathcal{L}_M \times \mathbb{R}$ by
$$r_*(v,\tau)=(r_*v,\tau), \qquad (v,\tau) \in \mathcal{L}_M \times \mathbb{R},
\quad r \in T^m.$$
\begin{lemma}
The delayed Rabinowitz action functional $\mathcal{A}_1$ is invariant under
the action of the torus $T^m$ on $\mathcal{L}_M \times \mathbb{R}$.
\end{lemma}
\textbf{Proof: } The area functional
$$\mathcal{L}_M \to \mathbb{R}, \quad v=(v_1,\ldots, v_m)
\to \int_{S^1} v^*\lambda=\sum_{i=1}^m \int_{S^1} v_i^*\lambda_i$$
is invariant under the action of $T^m$, since each
functional
$$\mathcal{L}_{M_i} \to \mathbb{R}, \quad v_i \mapsto \int_{S^1} v_i^* \lambda_i$$
is invariant under the circle action on $\mathcal{L}_{M_i}$. Moreover, each
averaged function
$$\overline{H}_i \colon \mathcal{L}_{M_i} \to \mathbb{R}$$
is invariant under the circle action as well, so that $f \circ \overline{H}$ is
again invariant under the action of $T^m$. This proves the lemma. \hfill $\square$
\\ \\
In contrast to the delayed Rabinowitz action functional the undelayed one
is in general only invariant under the standard circle action of
$\mathcal{L}_M \times \mathbb{R}$ given by reparametrising the loop which we obtain
from the torus action by embedding the circle diagonally into the torus
$$S^1 \to T^m, \quad r \mapsto (r,r,\ldots,r).$$
An exception is the case where the function $f$ is linear. Since the integral
is linear one has in this case $f \circ \overline{H}=\overline{f \circ H}$, so that
the delayed and the undelayed Rabinowitz action functionals coincide and are therefore
both invariant under the torus action. 
\\ \\
Critical points of $\mathcal{A}_1$ are solution of the Hamiltonian delay equation
\begin{equation}\label{critdel}
\left.\begin{array}{cc}
\partial_t v(t)=\tau X_{df(\overline{H}(v))H}(v(t)),& t \in S^1\\
f \circ \overline{H}(v)=0.
\end{array}\right\}
\end{equation}
This equation looks like a delay equation but it is actually an ODE in disguise and we will see that solutions of
(\ref{critdel}) are actually precisely the solution of (\ref{crit}). Even more is true
as Theorem\,A from the introduction says.
We interpolate between the two functional as follows. Namely for $r \in [0,1]$ we set
$$\mathcal{A}_r \colon \mathcal{L} \times \mathbb{R} \to \mathbb{R}, \quad
(v,\tau) \mapsto -\int_{S^1}v^*\lambda+\tau\Big(r f \circ \overline{H}(v)+(1-r)
\overline{f \circ H}(v)\Big).$$
We are now in position to prove Theorem\,A from the introduction.
\\ \\
\textbf{Proof of Theorem\,A: } Suppose that $(v,\tau) \in \mathrm{crit}(\mathcal{A}_r)$. Then
$(v,\tau)$ is a solution of the problem
\begin{equation}\label{critint}
\left.\begin{array}{cc}
\partial_t v(t)=\tau X_{rdf(\overline{H}(v))H+(1-r)f H}(v(t)),& t \in S^1\\
rf \circ \overline{H}(v)+(1-r) \overline{f \circ H}(v)=0.
\end{array}\right\}
\end{equation}
Abbreviate by
$$f_i:=\frac{\partial f}{\partial x_i}, \quad 1 \leq i \leq m$$
the partial derivatives of $f$. We rewrite the first equation in (\ref{critint})
componentwise as
\begin{eqnarray}\label{comp}
\partial_t v_i(t)
&=&\tau r f_i\Big(\overline{H}_1(v_1),\ldots \overline{H}_m(v_m)\Big)X_{H_i}(v_i(t))\\
\nonumber
& &+\tau(1-r)f_i\Big(H_1(v_1(t)),\ldots H_m(v_m(t))\Big)X_{H_i}(v_i(t)).
\end{eqnarray}
By preservation of energy it follows from (\ref{comp}) that $H_i(v_i)$ is constant
so that we have
\begin{equation}\label{mean}
H_i(v_i)(t)=\overline{H}_i(v_i), \quad t \in S^1.
\end{equation}
Plugging this into (\ref{comp}) we obtain
$$\partial_t v_i(t)
=\tau f_i\Big(H_1(v_1(t)),\ldots H_m(v_m(t))\Big)X_{H_i}(v_i(t))$$
which is independent of $r$. This implies that the first equation of (\ref{critint})
does not depend on the homotopy parameter. 
\\ \\
From (\ref{mean}) we infer further that
$$f \circ \overline{H}(v)=\overline{f \circ H}(v)$$
so that the second equation in (\ref{critint}) becomes
$$\overline{f \circ H}(v)=0$$
which is independent of $r$ as well. This shows that the critical set does not depend
on the homotopy parameter. 
\\ \\
It remains to explain why $\mathcal{A}_r$ is constant on the critical set.
If $(v,\tau) \in \mathrm{crit}(\mathcal{A}_r)$ it holds that
$$\mathcal{A}_r(v,\tau)=-\int_{S^1}v^* \lambda.$$
This expression does not depend on $r$ and the proof of the lemma is complete.
\hfill $\square$

\section{The delayed fundamental lemma}\label{delfund}

In this section we explain the restricted contact type condition for several
particles needed in Theorem\,B from the introduction. We then show that
under the restricted contact type condition the fundamental lemma in Rabinowitz
Floer homology continues to hold for the delayed case. Having the fundamental
lemma established the compactness proof for gradient flow lines is standard
and is briefly recalled at the end of the section.
\\ \\
Note that the primitive $\lambda$ of the symplectic form $\omega$ on
$M$, uniquely determines a Liouville vector field $Y=Y_\lambda$ which is implicitly defined by the condition
$$\lambda=\omega(Y,\cdot).$$
The first hypothesis we want to assume throughout this section is
\begin{description}
 \item[(H1)] The function $H_f$ has $0$ as a regular value and the energy hypersurface
$$\Sigma:=H_f^{-1}(0)$$
is compact and positively transverse to the Liouville vector field $Y$, in the sense
that
$$\lambda\big(X_{H_f}\big)\big|_{\Sigma}>0.$$
\end{description}
It follows from Hypothesis (H1) that  
the restriction of $\lambda$ to $\Sigma$ is a contact form on $\Sigma$ and at
every point in $\Sigma$ the Hamiltonian vector field $X_{H_f}$ is proportional
to the Reeb vector field by a positive proportionality constant. Our second hypothesis is 
\begin{description}
 \item[(H2)] The functions $H_i$ are constant outside of a compact set for
 $1 \leq i \leq m$. 
\end{description}
It follows from Hypothesis (H2) that the  Hamiltonian vector fields $X_{H_i}$ have compact support. Our third hypothesis is
\begin{description}
 \item[(H3)] The exact symplectic manifolds $(M_i,\lambda_i)$ are completions of
 Liouville domains for $1 \leq i \leq m$. 
\end{description}
In view of Hypothesis (H3) we choose an $\omega$-compatible almost complex structure on
$M$ which has the property that outside of a compact set
$$J=\bigoplus_{i=1}^m J_i$$
where $J_i$ is an $\omega_i$ compatible almost complex structure on
$M_i$ which is SFT-like outside a compact subset of $M_i$. 
On $\mathcal{L} \times \mathbb{R}$ we consider the
$L^2$-metric $g=g_J$ which at a point $(v,\tau) \in \mathcal{L} \times \mathbb{R}$ is given for
tangent vectors
$$(\hat{v}_1,\hat{\tau}_1), (\hat{v}_2,\hat{\tau}_2) \in T_v \mathcal{L} \times
\mathbb{R}$$
is given by 
$$g\big((\hat{v}_1,\hat{\tau}_1), (\hat{v}_2,\hat{\tau}_2)\big)
=\int_0^1 \omega \big(\hat{v}_1(t),J(v(t))\hat{v}_2(t)\big)+\hat{\tau}_1 \cdot
\hat{\tau}_2.$$
With respect to this metric the gradient $\nabla \mathcal{A}_r=\nabla_J \mathcal{A}_r$
at a point $(v,\tau) ´\in \mathcal{L} \times \mathbb{R}$ becomes
$$\nabla \mathcal{A}_r(v,\tau)=\left(\begin{array}{c}
J(v)\big(\partial_t v-\tau
X_{rdf(\overline{H}(v))H+(1-r)fH}(v)\big)\\
rf \circ
\overline{H}(v)+(1-r)\overline{f \circ H}(v)\end{array}\right)
\in T_v \mathcal{L} \times \mathbb{R}.$$
Hence gradient flow lines of $\nabla \mathcal{A}_r$ are solutions
$(v,\tau) \in C^\infty(\mathbb{R} \times S^1,M) \times C^\infty(\mathbb{R},\mathbb{R})$
of the problem
\begin{equation}\label{grad}
\left.\begin{array}{c}
\partial_s v+J(v)\big(\partial_t v-\tau
X_{rdf(\overline{H}(v))H+(1-r)fH}(v)\big)=0\\
\partial_s \tau+rf \circ
\overline{H}(v)+(1-r)\overline{f \circ H}(v)=0.
\end{array}\right\}
\end{equation}
In the following we denote by $||\nabla \mathcal{A}_r||$ the norm of the gradient 
with respect to the metric $g_J$. The next lemma tells us that if the norm of
the gradient is small we can bound the Lagrange multiplier $\tau$ in terms of
the action. With the help of this lemma the compactness proof is standard
and follows along the same lines as in \cite{cieliebak-frauenfelder1}. Therefore
in \cite{albers-frauenfelder} this lemma is referred to as the \emph{fundamental lemma in Rabinowitz Floer homology}. In the delayed case it requires quite some additional
work to establish it compared to the nondelayed one in \cite{cieliebak-frauenfelder1}.
\begin{lemma}\label{fundlem}
There exists a constant $c>0$ such that the following implication holds for
$(v,\tau) \in \mathcal{L} \times \mathbb{R}$ and $r \in [0,1]$
\begin{equation}\label{mainimpl}
||\nabla \mathcal{A}_r(v,\tau)|| \leq \frac{1}{c}\quad \Longrightarrow \quad
|\tau| \leq c\big(|\mathcal{A}_r(v,\tau)|+1\big).
\end{equation}
\end{lemma}
\textbf{Proof: } In view of Hypotheses (H1) and (H2) there exists a constant $\kappa>0$ and a constant
\begin{equation}\label{epsi}
0<\epsilon \leq \frac{\kappa}{4}
\end{equation}
such that
the following implication holds true
\begin{equation}\label{impl1}
|H_f(z)| \leq \epsilon \quad \Longrightarrow \quad
\lambda \big(X_{H_f}(z)\big) \geq \kappa.
\end{equation}
Since outside of a compact set the $\omega_i$-compatible almost complex structures
$J_i$ are
SFT-like the norm of $\lambda$ is uniformly bounded and hence there exists a constant
$L>0$ such that
\begin{equation}\label{lambound}
||\lambda_z|| \leq L, \quad z \in M.
\end{equation}
The functions $H_i \colon M_i \to \mathbb{R}$ we freely interpret as well
as functions on $M$ in the obvious way by pulling them back to $M$ under
the canonical projection $\pi_i \colon M \to M_i$, namely for $z=(z_1,\ldots,z_m) \in M$
we have
$$H_i(z):=H_i(z_i).$$
Since outside of a compact set the metric induced from the $\omega$-compatible
almost complex structure $J$ is of product type the norm $||dH_i(z)||$ there only
depends on $z_i$ and since $H_i$ on $M_i$ is constant outside of a compact set,
there exists a further constant $C$ such that
\begin{equation}\label{conC}
||dH_i(z)|| \leq C, \qquad z \in M, \quad 1 \leq i \leq m.
\end{equation}
In view of Hypothesis (H1) the image of $H$ in $\mathbb{R}^m$ is compact and 
therefore $f|_{\mathrm{im}(H)}$ is uniformly continuous. Therefore there exists
$\delta>0$ with the property that for $h_1,h_2 \in \mathrm{im}(H)$ one has the following
implication
\begin{equation}\label{glst}
||h_1-h_2||\leq \delta \quad \Longrightarrow \quad \big|f(h_1)-f(h_2)\big|\leq \frac{\epsilon}{3}.
\end{equation}
Maybe after shrinking $\delta$ we can in view of (\ref{impl1}) assume 
that the following implication holds
\begin{equation}\label{impl2}
|H_f(z)| \leq \epsilon, \,\,|h-H(z)| \leq \delta \quad
\Longrightarrow \quad \lambda \big(X_{df(h)H}(z)\big) \geq \frac{\kappa}{2}.
\end{equation}
If $v \in \mathcal{L}$ we define its $H$-oscillation
$$\mathfrak{o}(v):=\max\Big\{\big|\big|H(v(t_1))-H(v(t_2))\big|\big|: t_1,t_2 \in S^1\Big\}.$$
We prove now the lemma in three steps.
\\ \\
\textbf{Step\,1:} \emph{Assume that $\mathfrak{o}(v) \leq \delta$, then the following
implication holds}
\begin{equation}\label{est0}
\big|\big|\nabla \mathcal{A}_r(v,\tau)\big|\big| \leq \frac{2\epsilon}{3}
\quad \Longrightarrow \quad
|\tau| \leq \frac{6}{\kappa}\bigg(\big|\mathcal{A}_r(v,\tau)\big|+
\frac{2\epsilon L}{3 }\bigg).
\end{equation}
For $t \in S^1$ we estimate
\begin{eqnarray}\label{durch1}
\bigg|\bigg|H(v(t))-\overline{H}(v)\bigg|\bigg|&=&
\bigg|\bigg|H(v(t))-\int_0^1 H(v(s))ds\bigg|\bigg|\\ \nonumber
&=&\bigg|\bigg|\int_0^1\Big(H(v(t))-H(v(s))\Big)ds\bigg|\bigg|\\ \nonumber
&\leq&\int_0^1 \big|\big| H(v(t))-H(v(s))\big|\big|ds\\ \nonumber
&\leq&\int_0^1 \mathfrak{o}(v)ds\\ \nonumber
&=&\mathfrak{o}(v)\\ \nonumber
&\leq& \delta.
\end{eqnarray}
In view of (\ref{glst}) this implies
\begin{equation}\label{durchi}
\big|f\big(H(v(t))\big)-f\big(\overline{H}(v)\big)\big| \leq \frac{\epsilon}{3}.
\end{equation}
Since $t \in S^1$ is arbitrary we deduce from that
\begin{eqnarray}\label{durch2}
\bigg|\overline{f \circ H}(v)-f\circ \overline{H}(v)\bigg|
&=&\bigg|\int_0^1 f\big(H(v(t))\big)dt-f\big(\overline{H}(v)\big)\bigg|\\ \nonumber
&\leq& \int_0^1 \big|f\big(H(v(t))\big)-f\big(\overline{H}(v)\big)\big|dt\\ \nonumber
&\leq&\frac{\epsilon}{3}.
\end{eqnarray}
We now subdivide the proof of Step\,1 into three substeps. 
\\ \\
\textbf{Step\,1a: } \emph{Assume that 
$|H_f(v(t))| \leq \epsilon$ for every $t \in S^1$, then}
\begin{equation}\label{est1}
|\tau| \leq \frac{6}{\kappa}\bigg(\big|\mathcal{A}_r(v,\tau)\big|+
L\big|\big|\nabla \mathcal{A}_r(v,\tau)\big|\big|\bigg).
\end{equation}
Suppose that $t \in S^1$. In view of the assumption of Step\,1a it follows from
(\ref{impl1}) that
$$\lambda\big(X_{H_f}(v(t)\big)\geq \kappa.$$
Using additionally (\ref{durch1}) we infer with (\ref{impl2}) that
$$\lambda\big(X_{df(\overline{H}(v))H}(v(t))\big)\geq \frac{\kappa}{2}.$$
Therefore we obtain
\begin{eqnarray}\label{reeb}
& &\lambda\big(X_{rdf(\overline{H}(v))H+(1-r)fH}(v(t))\big)
\\ \nonumber
&=&r \lambda\big(X_{df(\overline{H}(v))H}(v(t))\big)+(1-r)\lambda\big(X_{H_f}(v(t)\big)
\\ \nonumber
&\geq& \frac{r\kappa}{2}+(1-r)\kappa\\ \nonumber
&\geq& \frac{\kappa}{2}.
\end{eqnarray}
In view of the assumption of Step\,1a we infer that
$$\big|\overline{f \circ H}(v)\big| \leq \epsilon.$$
Combining this estimate with (\ref{durch2}) we get
\begin{eqnarray*}
\big|f \circ \overline{H}(v)\big|&\leq&
\big|f \circ \overline{H}(v)-\overline{f \circ H}(v)\big|+
\big|\overline{f \circ H}(v)\big|\\
&\leq& \frac{\epsilon}{3}+\epsilon\\
&=&\frac{4\epsilon}{3}.
\end{eqnarray*}
From the above two estimates we infer
\begin{eqnarray}\label{mit}
\big|r f \circ \overline{H}(v)+(1-r)\overline{f \circ H}(v)\big|
&\leq& r\big|f \circ \overline{H}(v)\big|+(1-r)\big|\overline{f \circ H}(v)\big|\\
\nonumber
&\leq& \frac{4\epsilon r}{3}+(1-r)\epsilon\\ \nonumber
&\leq& \frac{4\epsilon}{3}.
\end{eqnarray}
Using the inequalities (\ref{reeb}) and (\ref{mit}), the uniform bound
$L$ on the one-form $\lambda$ from (\ref{lambound}) and remembering that in
(\ref{epsi}) we have chosen $\epsilon \leq \tfrac{\kappa}{4}$ we estimate
\begin{eqnarray*}
\big|\mathcal{A}_r(v,\tau)\big|&\geq&\bigg|\int_{S^1}v^*\lambda\bigg|-
\big|\tau\big|\cdot \big|rf \circ \overline{H}(v)+(1-r)\overline{f \circ H}(v)\big|\\
&\geq&\bigg|\int_0^1\lambda\Big(\tau X_{r df(\overline{H}(v))H+
(1-r)fH}(v)\Big)dt\bigg|\\
& &-\bigg|\int_0^1 \lambda\Big(\partial_t(v)-\tau X_{r df(\overline{H}(v))H+
(1-r)fH}(v)\Big)dt\bigg|-\frac{4 \epsilon |\tau|}{3}\\
&\geq&\frac{\kappa |\tau|}{2}-L\bigg|\bigg|\partial_t v-\tau X_{r df(\overline{H}(v))H+
(1-r)fH}(v)\bigg|\bigg|_{L^1(S^1)}-\frac{\kappa |\tau|}{3}\\
&\geq& \frac{\kappa |\tau|}{6}-L\bigg|\bigg|\partial_t v-\tau X_{r df(\overline{H}(v))H+
(1-r)fH}(v)\bigg|\bigg|_{L^2(S^1)}\\
&\geq&\frac{\kappa |\tau|}{6}-L \big|\big|\nabla \mathcal{A}_r(v,\tau)\big|\big|
\end{eqnarray*}
from which (\ref{est1}) follows. This proves Step\,1a.
\\ \\
\textbf{Step\,1b: } \emph{Assume that there exists $t_0 \in S^1$ such that
$|H_f(v(t_0)|>\epsilon$, then}
\begin{equation}\label{est2}
\big|\big|\nabla \mathcal{A}_r(v,\tau)\big|\big|> \frac{2\epsilon}{3}.
\end{equation}
Suppose that $t \in S^1$. Since the $H$-oscillation is bounded by $\delta$ we
have
$$||H(v(t))-H(v(t_0)|| \leq \delta$$
and therefore we estimate with the help of (\ref{glst})
\begin{eqnarray*}
\big|H_f(v(t)\big|&=&\big|f(H(v(t)))\big|\\
&\geq&\big|f(H(v(t_0)))\big|-\big|f(H(v(t)))-f(H(v(t_0)))\big|\\
&>&\epsilon-\frac{\epsilon}{3}\\
&=&\frac{2\epsilon}{3}.
\end{eqnarray*}
Since the circle is connected we either have
\begin{equation}\label{case1}
H_f(v(t)) > \frac{2\epsilon}{3}, \qquad \forall\,\,t \in S^1
\end{equation}
or 
\begin{equation}\label{case2}
H_f(v(t)) < -\frac{2\epsilon}{3}, \qquad \forall\,\,t \in S^1.
\end{equation}
We first discuss (\ref{case1}). In this case we have
\begin{equation}\label{a1}
\overline{f \circ H}(v) > \frac{2\epsilon}{3}.
\end{equation}
Moreover, it holds that
$$H_f(v(t_0))> \epsilon$$
and hence applying (\ref{durchi}) for $t=t_0$ we infer that
\begin{equation}\label{a2}
f \circ \overline{H}(v) > \frac{2\epsilon}{3}.
\end{equation}
From (\ref{a1}) and (\ref{a2}) we infer that
$$r f\circ \overline{H}(v)+(1-r)\overline{f \circ H} > \frac{2\epsilon}{3}$$
from which (\ref{est2}) follows. 
This proves Step\,1b in case (\ref{case1}) holds. The case (\ref{case2}) is similar.
There we infer
$$\overline{f \circ H}(v)<-\frac{2\epsilon}{3},\qquad
f\circ \overline{H}(v)<-\frac{2\epsilon}{3}$$
from which follows
$$r f\circ \overline{H}(v)+(1-r)\overline{f \circ H} < -\frac{2\epsilon}{3}$$
which again implies (\ref{est2}). This finishes the proof of Step\,1b. 
\\ \\
\textbf{Step\,1c:} \emph{We prove Step\,1.}
\\ \\
We assume that
$$\big|\big|\nabla \mathcal{A}_r(v,\tau)\big|\big| \leq \frac{2\epsilon}{3}.$$
In view of Step\,1b this implies that $|H_f(v(t))| \leq \epsilon$ for
every $t \in S^1$. Therefore we can apply Step\,1a and infer that
\begin{eqnarray*}
|\tau| &\leq&\frac{6}{\kappa}\bigg(\big|\mathcal{A}_r(v,\tau)\big|+
L\big|\big|\nabla \mathcal{A}_r(v,\tau)\big|\big|\bigg)\\
&\leq&\frac{6}{\kappa}\bigg(\big|\mathcal{A}_r(v,\tau)\big|+
\frac{2\epsilon L}{3}\bigg).
\end{eqnarray*}
This finishes the proof of Step\,1. 
\\ \\
In the second step we treat the case of large $H$-oscillation. For that 
recall the constant $C$ which appeared in (\ref{conC}).
\\ \\
\textbf{Step\,2: } \emph{Assume that $\mathfrak{o}(v)>\delta$, then}
\begin{equation}\label{est3}
\big|\big|\nabla \mathcal{A}_r(v,\tau)\big|\big|>\frac{\delta}{C\sqrt{m}}.
\end{equation}
By definition of the oscillation there exist times $t_0, t_1 \in S^1$ such that
$$\big|\big|H(v(t_1))-H(v(t_0))\big|\big|>\delta.$$
In particular, there exists $1 \leq i \leq m$ such that
$$\big|H_i(v_i(t_1))-H_i(v_i(t_0))\big|> \frac{\delta}{\sqrt{m}}.$$
Maybe after interchanging the roles of $t_1$ and $t_0$ we can even assume that
$$H_i(v_i(t_1))-H_i(v_i(t_0))> \frac{\delta}{\sqrt{m}}.$$
By going to the universal cover $\mathbb{R}$ of the circle $S^1=\mathbb{R}/\mathbb{Z}$
we interpret $t_0$ and $t_1$ as real numbers satisfying
$$t_0<t_1<t_0+1.$$
Abbreviate
$$P:=\big\{t \in [t_0,t_1]: dH_i(v_i(t))\partial_t v_i(t) \geq 0\big\}.$$
We estimate
\begin{eqnarray*}
& &\frac{\delta}{\sqrt{m}}\\
&<&H_i(v_i(t_1))-H_i(v_i(t_0))\\
&=&\int_{t_0}^{t_1} dH_i(v_i(t))\partial_t v_i(t)dt\\
&\leq&\int_P dH_i(v_i(t))\partial_t v_i(t)dt\\
&=&\int_PdH_i(v_i(t))\bigg(\partial_t v_i(t)-\tau \Big(r f_i
(\overline{H}(v))+(1-r)f_i(H(v(t))\Big)X_{H_i}(v_i(t))\bigg)dt\\
&=&\int_P dH_i(v)\Big(\partial_t v-\tau X_{rdf(\overline{H}(v))H+(1-r)fH}(v)\Big)dt\\
&\leq&C \int_P \big|\big|\partial_t v-\tau X_{rdf(\overline{H}(v))H+(1-r)fH}(v)
\big|\big|dt\\
&\leq& C\big|\big|\partial_t v-\tau X_{rdf(\overline{H}(v))H+(1-r)fH}(v)
\big|\big|_{L^1(S^1)}\\
&\leq& C\big|\big|\partial_t v-\tau X_{rdf(\overline{H}(v))H+(1-r)fH}(v)
\big|\big|_{L^2(S^1)}\\
&\leq&C\big|\big|\nabla \mathcal{A}_r(v,\tau)\big|\big|.
\end{eqnarray*}
This implies (\ref{est3}) and finishes the proof of Step\,2. 
\\ \\
\textbf{Step\,3: } \emph{We prove the lemma.}
\\ \\
We define
$$c:=\max\bigg\{\frac{C\sqrt{m}}{\delta},\frac{3}{2\epsilon},\frac{6}{\kappa},
\frac{4\epsilon L}{\kappa}\bigg\}$$
and show that with this choice of the constant $c$ the implication
(\ref{mainimpl}) holds true. For this purpose suppose that
\begin{equation}\label{hiest}
\big|\big|\nabla \mathcal{A}_r(v,\tau)\big|\big| \leq \frac{1}{c}.
\end{equation}
This implies that
$$\big|\big|\nabla \mathcal{A}_r(v,\tau)\big|\big| \leq \frac{\delta}{C \sqrt{m}}.$$ 
From Step\,2 we infer that
$$\mathfrak{o}(v) \leq \delta.$$
Therefore we deduce from Step\,1 that the implication (\ref{est0}) holds true. From
(\ref{hiest}) we further deduce that
$$\big|\big|\nabla \mathcal{A}_r(v,\tau)\big|\big| \leq \frac{2\epsilon}{3}.$$
Hence we obtain from (\ref{est0}) that
\begin{eqnarray*}
\big|\tau\big| &\leq&\frac{6}{\kappa}\bigg(\big|\mathcal{A}_r(v,\tau)\big|+
\frac{2\epsilon L}{3}\bigg)\\
&=&\frac{6}{\kappa}\big|\mathcal{A}_r(v,\tau)\big|+\frac{4\epsilon L}{\kappa}\\
&\leq&c\big|\mathcal{A}_r(v,\tau)\big|+c\\
&=&c\big(\big|\mathcal{A}_r(v,\tau)\big|+1\big).
\end{eqnarray*}
This proves the implication (\ref{mainimpl}) and the lemma follows. \hfill $\square$
\\ \\
Having Lemma~\ref{fundlem} at our disposal the compactness proof for the gradient
flow lines now follows precisely the same scheme as the compactness proof in
\cite{cieliebak-frauenfelder1}, which then proves Theorem\,B from the introduction.
We sketch the main steps.
\\ \\
\textbf{Sketch of proof of Theorem\,B: } The fundamental lemma allows to bound
the Lagrange multiplier in terms of action. Since the action stays in the bounded
interval $[a,b]$ these leads to a uniform bound on the Lagrange multiplier. Because
outside of a compact set the almost complex structures $J_i$ are SFT-like and
the Hamiltonian vector fields $X_{H_i}$ vanish, the maximum principle tells us
that gradient flow lines have to stay in a compact subset of $M$. Moreover, there
derivatives cannot explode since the symplectic form $\omega$ is exact and therefore
there is no bubbling. With these uniform bounds the compactness theorem follows
from elliptic regularity of the Cauchy-Riemann operator. \hfill $\square$


\begin{thebibliography}{99}

\bibitem{albers-frauenfelder} P.\,Albers, U.\,Frauenfelder, \emph{Rabinowitz
Floer homology: a survey.} Global differential geometry, 437--461, Springer
Proc.\,Math., \textbf{17}, Springer, Heidelberg (2012).
\bibitem{albers-frauenfelder-schlenk} P.\,Albers, U.\,Frauenfelder, F.\,Schlenk,
\emph{Hamiltonian delay equations -- examples and  lower bound for the number
of periodic solutions}, Adv.\,Math. \textbf{373}, 107319 (2020). 
\bibitem{cieliebak-frauenfelder1} K.\,Cieliebak, U.\,Frauenfelder, \emph{A Floer homology for
exact contact embeddings}, Pacific J.\,Math. \textbf{239} (2009), no.\,2, 251--316.
\bibitem{cieliebak-frauenfelder2} K.\,Cieliebak, U.\,Frauenfelder,
\emph{Spectral numbers in Tate Rabinowitz Floer homology}, preprint.
\bibitem{cieliebak-frauenfelder-volkov} K.\,Cieliebak, U.\,Frauenfelder, E.\,Volkov,
\emph{A variational approach to frozen planet orbits in Helium}, to appear in
Ann.\,Inst.\,H.\,Poincar\'e.
\bibitem{frauenfelder0} U.\,Frauenfelder, \emph{Helium and Hamiltonian delay equations},
Israel J.\,Math. \textbf{246}, no.\,1, 239--260 (2021). 
\bibitem{frauenfelder} U.\,Frauenfelder, \emph{The Stark problem as a concave
toric domain}, Geom.\,Dedicata \textbf{217}, no.\,1, Paper No.\,10 (2023).
\bibitem{gutzwiller} M.\,Gutzwiller, \emph{Chaos in classical and quantum mechanics},
Springer-Verlag, New York (1990).
\bibitem{mohebbi} A.\,Mohebbi, \emph{The special concave toric domain for the
rotating Kepler problem}, arXiv:2108.04581
\bibitem{pinzari} G.\,Pinzari, \emph{Proof of a conjecture by H.\,Dullin and
R.\,Montgomery}, arXiv:2209.07097
\bibitem{ramos} V.\,Ramos, \emph{Symplectic embeddings and the Lagrangian bidisk},
Duke Math.\,J. \textbf{166}, no.\,9, 1703--1738 (2017). 

\end{thebibliography}
\end{document}